\documentclass[12pt]{article}
\usepackage{amssymb,amsthm,amsmath,latexsym}
\usepackage{url}
\usepackage{lscape}
\newtheorem{thm}{Theorem}[section]

\theoremstyle{remark}

\title{There exists a partitioned balanced tournament design of side nine.}

\author{%
Makoto Araya\thanks{Department of Computer Science, Shizuoka University}\quad and \quad Naoya Tokihisa\thanks{Graduate School of Integrated Science and Technology, Shizuoka University}
}
\date{May 17, 2018}
\begin{document}

\maketitle

\begin{abstract}
We construct a partitioned balanced tournament design of side nine.
\end{abstract}

\section{Introduction}

A {\it partitioned balanced tournament design of side $n$}, \textrm{PBTD($n$)}, 
defined on a $2n$-set $V$ is an arrangement of the ${2n \choose n}$
distinct unordered pairs of the elements of $V$ into an $n\times (2n-1)$ arrays such that
\begin{enumerate}
\item every element of $V$ is contained in precisely one cell of each column,
\item every element of $V$ is contained in at most two cells of any row,
\item each row contains a factor in the first $n$ columns, and
\item each row contains a factor in the last $n$ columns.
\end{enumerate}

E.R. Lamken prove the following theorem.

\begin{thm}[\cite{hb}]
There exists a \textrm{PBTD($n$)} for $n$ a positive integer, $n \ge 5$, except possibly for $n \in \{9,11,15\}$.
\end{thm}

By a computer calculation, we construct a PBTD(9) in Table~$1$.
Hence the following theorem holds.
\begin{thm}
There exists a partitioned balanced tournament design of side nine.
\end{thm}

\begin{landscape}
\begin{table}
$\begin{array}{|c|c|c|c|c|c|c|c|c|c|c|c|c|c|c|c|c|}
\hline
2, 16  &  3, 17  &  4, 6  &  5, 7  &  8, 10  &  9, 11  &  12, 14  &  13, 15  &  0, 1 & 2, 5  &  3, 4  &  6, 15  &  7, 14  &  8, 11  &  9, 10  &  12, 16  &  13, 17\\\hline
0, 4  &  1, 5  &  7, 9  &  6, 8  &  11, 13  &  10, 12  &  15, 17  &  14, 16  &  2, 3 & 0, 16  &  1, 17  &  4, 8  &  5, 9  &  6, 13  &  7, 12  &  10, 15  &  11, 14\\ \hline
1, 3  &  0, 2  &  10, 13  &  11, 12  &  14, 17  &  15, 16  &  6, 9  &  7, 8  &  4, 5 & 6, 10  &  7, 11  &  1, 16  &  0, 17  &  9, 12  &  8, 13  &  2, 14  &  3, 15\\ \hline
10, 14  &  11, 15  &  0, 8  &  1, 9  &  2, 4  &  3, 5  &  13, 16  &  12, 17  &  6, 7 & 3, 13  &  2, 12  &  9, 17  &  8, 16  &  4, 14  &  5, 15  &  0, 11  &  1, 10\\ \hline
5, 6  &  4, 7  &  2, 17  &  3, 16  &  12, 15  &  13, 14  &  0, 10  &  1, 11  &  8, 9 & 4, 11  &  5, 10  &  2, 13  &  3, 12  &  0, 15  &  1, 14  &  7, 17  &  6, 16\\ \hline
8, 12  &  9, 13  &  1, 15  &  0, 14  &  5, 16  &  4, 17  &  3, 7  &  2, 6  &  10, 11 & 1, 12  &  0, 13  &  5, 14  &  4, 15  &  7, 16  &  6, 17  &  3, 8  &  2, 9\\\hline
9, 15  &  8, 14  &  11, 16  &  10, 17  &  3, 6  &  2, 7  &  1, 4  & 0, 5  &  12, 13 &  9, 14  &  8, 15  &  3, 11  &  2, 10  &  5, 17  &  4, 16  &  1, 6  &  0, 7\\ \hline
11, 17  &  10, 16  &  5, 12  &  4, 13  &  1, 7  &  0, 6  &  2, 8  &  3, 9  &  14, 15 & 8, 17  &  9, 16  &  7, 10  &  6, 11  &  1, 2  &  0, 3  &  5, 13  &  4, 12\\ \hline
7, 13  &  6, 12  &  3, 14  &  2, 15  &  0, 9  &  1, 8  &  5, 11  &  4, 10  &  16, 17 & 7, 15  &  6, 14  &  0, 12  &  1, 13  &  3, 10  &  2, 11  &  4, 9  &  5, 8\\\hline
\end{array}$
\caption{a partitioned balanced tournament design of side nine}
\end{table}
\end{landscape}

\end{document}